\documentclass[12pt]{amsart}
\usepackage{amsmath}
\usepackage{enumerate}
\usepackage{amssymb}
\usepackage{amsbsy}
\usepackage{amsfonts}
\theoremstyle{plain}
\newtheorem{thm}{Theorem}[section]
\newtheorem{prop}[thm]{Proposition}

\newtheorem{cor}[thm]{Corollary}
\theoremstyle{definition}

\newtheorem{rem}[thm]{Remark}

\numberwithin{equation}{section}

\newcommand{\sm}{\begin{smallmatrix}}
\newcommand{\esm}{\end{smallmatrix}}

\newcommand{\C}{\mathbb{C}}

\newfont{\FieldFont}{msbm10 scaled\magstep1}

\newcommand{\pf}{\noindent\bf Proof }

\def\C{\mathbb C}

\begin{document}

\title[Non vanishing of Central values of modular L-functions ]{
Non vanishing of Central values of modular L-functions for Hecke eigenforms of level one}

\author{D. Choi }

\address{School of Liberal Arts and Sciences, Korea Aerospace University, Goyang, Gyeonggi 412-791, Korea}
\email{choija@kau.ac.kr}

\author{Y. Choie}

\address{Dept of Mathematics and  PMI, POSTECH, Pohang , Korea 790-784}
\email{yjc@postech.ac.kr}


\begin{abstract}
Let $F(z)=\sum_{n=1}^{\infty}a(n)q^n$ be a newform of weight $2k$
and level $N$  with a trivial character, and assume that
$F(z)$ is a non-zero eigenform of all Hecke operators. For $x>0,$
let
 $$\mathcal{N}_F(x):=|\{
 D \, \mbox{\, \, fundamental } | \, |D|<x, (D,N)=1, L(F, D,k)\neq
 0\}|.$$
 A based on the Goldfeld's conjecture  one expects to have
$ \mathcal{N}_F(x)\gg x \,\, \, \, (x \rightarrow \infty). $
Kohnen \cite{Kohnen3} showed that if $k\geq 6$ is a even integer,
then for $x\gg 0$ there is a normalized Hecke eigenform $F$ of
level $1$ and weight $2k$ with the property that
$$\mathcal{N}_F(x)\gg_k x \,\, \, \, (x \rightarrow \infty).$$
In this paper, we extend the result in \cite{Kohnen3} to the case
when $k$ is any integer, in particular when $k$ is odd. So, we obtain that, when the level is 1, for each
integer $2k$ such that the dimension of cusp forms of weight $2k$
is not zero, there is a normalized Hecke eigenform $F$ of weight
$2k$ satisfying $ \mathcal{N}_F(x)\gg x \,\, \, \, (x \rightarrow
\infty).  $
\end{abstract}

 \maketitle

\today

\section{\bf{Introduction and statement of result}}

Let $F(z)=\sum_{n=1}^{\infty}a(n)q^n \in S_{2k}(N, \chi_0)$ be a
newform of weight $2k$ and level $N$  with a trivial character
$\chi_0$, and assume that $F(z)$ is a non-zero eigenform of all
Hecke operators.  For a fundamental discriminant $D,$ that is
$D=1$ or the discriminant of a quadratic field, we define the
$L-$function $L(F,D,s)$ of $F$ twisted with the quadratic
character $(\frac{D}{\bullet})$ by
$$\sum_{n\geq 1}\left(\frac{D}{n}\right)\frac{a(n)}{n^{s}}.$$
In this paper, we consider the central values of $L-$functions
$L(F,D,k)$.

It is well-known by Waldspurger \cite{wal}  that the central
critical values $L(F,D,k)$ are essentially proportional to the
squares of Fourier coefficients of the modular form of weight
$k+\frac{1}{2}$ corresponding to $F$ under Shimura correspondence.
\medskip

 On the other hand, for $x>0,$ consider the set
 $$\mathcal{N}_F(x):=|\{
 D \, \mbox{\, \, fundamental } | \, |D|<x, (D,N)=1, L(F, D,k)\neq
 0\}|.$$
A based on the Goldfeld's conjecture  one expected to have
\begin{equation}\label{conjecture}
\mathcal{N}_F(x)\gg x \,\, \, \, (x \rightarrow \infty)
\end{equation}

( In \cite{Gf} Goldfeld conjectured that
$$\sum_{D \text{\, \, fundamental, } \, |D|<x, \gcd{(D,N)}=1}
ord_{s=k}L(F,D,s)$$
$$\approx \frac{1}{2} |\{ D \, \text{\, fundamental } \, | \,
|D|<x, \gcd{(D,N)}=1 \} |\;  \, (x \rightarrow \infty).)$$

 In \cite{O-S} using the theory of Galois representations of
modular forms together with the results of  Friedberg and
Hoffstein \cite{F-H}, the following was proved
$$\mathcal{N}_F(x)\gg_F \frac{x}{\log{x}} \, \, \, (x\rightarrow
\infty).$$

In \cite{James} K. James  gave an example of a form $F$, by
looking at the difference of two special ternary theta series of
level $56$, satisfying $\mathcal{N}_F(x)\gg_F x$ $(x\rightarrow
\infty).$ Furthermore, Kohnen \cite{Kohnen3} showed that if $k\geq
6$ is an even integer, then there is a normalized Hecke eigenform
$F$ of level $1$ and weight $2k$ with the property that
$\mathcal{N}_F(x)\gg_k x$ for $x\gg 0$; in particular, it was
shown that $\mathcal{N}_{\Delta}(x)\gg x \, \, (x\rightarrow
\infty)$ where $\Delta$ is Ramanujan's $\Delta$-function of weight
$12.$  More precisely,  Kohnen  proved that if $N_{k,\Gamma_{1}}^+
(x)$ is the number of fundamental discriminant $ D,  \,  0<D<x, $
such that there exists a normalized Hecke eigenform $F \in S_{2k}
$ satisfying $L(F,D,k) \neq 0$, then for any positive $\epsilon$
$$\mathcal{N}_{k,\Gamma_1}^+(x) \geq \left(\frac{9}{16
\pi^2}-\epsilon\right)x \, \, \, \,  (x \gg_{\epsilon}0) .$$ Here,
$S_{2k}$ is the space of cusp forms of weight $2k$ on
$\Gamma_1=SL(2, \mathbb{Z}).$ In
this paper, we extend the result \cite{Kohnen3} to the case
when $k$ is any integer, in particular when $k$ is odd. To obtain these results,
we refine the argument of Kohnen by using an isomorphism \cite{Kohnen2}
from the spaces of modular forms of integral weight to the Kohnen
plus space of half integral weight modular forms.

\medskip

Let
 $N_{k,\Gamma_1}  (x)$ denote  the number of fundamental discriminants $
D,  \, |D|<x, $ such that there exists a normalized Hecke
eigenform $F \in S_{2k} $ satisfying $L(F,D,k) \neq 0$, then
  we have the following:

\begin{thm}\label{main} Suppose that $k\geq 9$ is odd.
Then for any positive $\epsilon$
\begin{equation*}
\mathcal{N}_{k,\Gamma_1}(x)\geq
\left(\frac{9}{16\pi^2}-\epsilon\right)x \, \, \, (x\gg_{\epsilon}
0).
\end{equation*}
\end{thm}

Let $d_k:=\dim(S_{2k} )$. With the result of Kohnen \cite{Kohnen3}
in which covers the case of $k \equiv 0\pmod{2},$ we state  the
following result:
\begin{thm}\label{main2}
Suppose that  $k$ is a positive integer such that $d_k\geq1$. Then
for any $\epsilon>0$
\begin{equation*}
\mathcal{N}_{k,\Gamma_1}(x)\geq \left(\frac{9}{16\pi^2}-\epsilon\right)x \,
\, \,  (x\gg_{\epsilon} 0).
\end{equation*}
\end{thm}

 Theorem \ref{main2} immediately implies that for each integer $k$ such that $d_k>0$ there is
a Hecke eigenform $F\in S_{2k}$ satisfying (\ref{conjecture}).

\begin{cor}\label{sp}
Suppose that $k$ is an integer such that $d_k\geq1$. Then there
exists a Hecke eigenform $F(z) \in S_{2k} $ such that for any
positive $\epsilon$
\begin{equation*}
\mathcal{N}_{F}(x)\geq \frac{1}{d_k}
\left(\frac{9}{16\pi^2}-\epsilon \right)x \, \, \, \, (x
\gg_{\epsilon } 0).
\end{equation*}
\end{cor}

Recently, Farmer and James \cite{F-J} proved that the
characteristic polynomial of the Hecke operator $T_2$ on $S_{2k}$
is irreducible over $\mathbb{Q}$ for $k\leq 1, 000$. If $K$ is the
field obtained from $\mathbb{Q}$ by adjoining the eigenvalues of
$T_2$, then the Galois group $G=Gal(K/\mathbb{Q})$ operates
transitively on the set of normalized Hecke eigenforms in
$S_{2k}$.  Using Theorem \ref{main} and the known fact
\cite{Shi2} that
$L(F^{\sigma},D,k)_{alg}=L(F,D,k)^{\sigma}_{alg}$ for all $\sigma
\in G$,  we have that every Hecke eigenform $F\in S_{2k}$
satisfies (\ref{conjecture}) for each integers $k$, $6\leq
k\leq1,000$. Here, $``alg"$ means ``algebraic part".

\begin{thm}
Suppose that $6\leq k \leq 1,000$ is an integer. Then every
normalized Hecke eigenform $F$ in $S_{2k}$ satisfies
\begin{equation}\label{1.5}
\mathcal{N}_{F}(x)\geq  \left(\frac{9}{16\pi^2}-\epsilon \right)x \, \, \, \,
(x \gg_{\epsilon } 0).
\end{equation}
\end{thm}

\begin{rem}
Maeda (\cite{Hi-Ma} Conjecture 1.2) made a conjecture that the
Hecke algebra of $S_{2k}$ over $\mathbb{Q}$ is simple, and that
its Galois closure over $\mathbb{Q}$ has Galois group $G$ the full
symmetric group. The conjecture implies that there is a single
Galois orbit of Hecke eigenforms in $S_{2k}$. Thus,  Maeda's
conjecture implies that every normalized Hecke eigenform $F$ in
$S_{2k}$ satisfies (\ref{1.5}).
\end{rem}

\bigskip

\section{\bf{Preliminaries}}
Let $q:=e^{2 \pi i z}$, where $z$ is in the complex upper half plane
$\mathbb{H}$. For an integer $k\geq2$ recall the normalized Eisenstein series
$E_{2k}(z):=1-\frac{4k}{B_{2k}}\sum_{n\geq 1} \sigma_{k-1}(n)q^n$  of
weight $2k$, and for an integer $r\geq1$  let
$H_{r+\frac{1}{2}}(z)=\sum_{N=0}^{\infty}H(r,N)q^N$ be
the Cohen-Eisenstein series of weight $r+\frac{1}{2}$ on
$\Gamma_0(4) $ (see \cite{Co}). Here, for each positive integer $N,$ define
$$h(r,N)=\left \{ \begin{array}{ll}
(-1)^{[\frac{r}{2}]}(r-1)!N^{r-\frac{1}{2}}2^{1-r}\pi^{-r}L(r,
\chi_{(-1)^rN}) & \mbox{if $(-1)^rN\equiv0$ or $1 \pmod{4} $},\\
0 & \mbox{if  $(-1)^rN \equiv 2 $ or $3\pmod{4}$},
\end{array}\right.$$
where  $\chi_D$ denotes the character $\chi_{D}(d)=(\frac{D}{d}).$
Furthermore, for $N\geq 1,$ define
$$H(r,N)=\left\{\begin{array}{lll} \sum_{d^2|N}
h(r,\frac{N}{d^2}) & \mbox{if $(-1)^r N \equiv 0$ or $1 \pmod{4}$},\\
\zeta(1-2r) & \mbox{if $N=0$},\\
0 & \mbox{otherwise}.\end{array} \right. $$

\medskip

The followings are proved in \cite{Co}:
\begin{prop}\label{eis}
\begin{enumerate}[(1)]
\item For $r\geq2$, Cohen-Eisenstein series $H_{r+\frac{1}{2}}(z)$ is a modular
form of weight $r+\frac{1}{2}$ on $\Gamma_0(4)$ and it is in
Kohnen plus condition, that is,
$$ H(r,N)=0 \mbox{\, if $(-1)^r N \neq 0,1 \pmod{4}$}. $$

\item  Let $a$ and $b$ be integers with $a \geq 1$. Suppose that $-b$ is a quadratic non residue of
$a$. Then the function
$$G_{a,b}(z):=\sum_{N\equiv b \pmod{a}}H(1,n)q^n$$
is a modular form of weight $\frac{3}{2}$ and character a over
$\Gamma_0(A)$, where we can take $A=4a^2$, and furthermore $A=a^2$
if $a$ is even.
\end{enumerate}
\end{prop}
For a nonnegative integer $k$ denote
$M_{k+\frac{1}{2}}(\Gamma_0(4))$ as the usual complex vector space
of cusp forms of weight $k+\frac{1}{2}$ on $\Gamma_0(4)$ with the
  trivial character. In \cite{Kohnen2}, Kohnen introduces
the plus space $M^{+}_{k+\frac{1}{2}}(\Gamma_0(4))$ of modular
forms $g(z)$ of weight $k+\frac{1}{2}$ on $\Gamma_0(4)$ with
a Fourier expansion of the form
\begin{equation*}g(z)=\sum_{(-1)^{k}n\equiv0,1
\pmod{4}} c(n)q^n.
\end{equation*}
and proved the following isomorphism (Proposition 1 in
\cite{Kohnen2}).

\begin{prop}\label{isomorphism}
Let $M_k$ be the space of modular forms of weight $k$ on
$\Gamma_1$. If $k$ is even, then the spaces $M_k\bigoplus M_{k-2}$
and $M_{k+\frac{1}{2}}^+(\Gamma_0(4))$ are isomorphic under the map
$$(f(z),h(z))\mapsto f(4z)\theta(z)+h(4z)H_{\frac{5}{2}}(z),$$
where $\theta(z)=1+2\sum_{n=1}^{\infty}q^{n^2}$.  If $k$ is odd,
then the spaces $M_{k-3}\bigoplus M_{k-5}$ and
$M_{k+\frac{1}{2}}^+(\Gamma_0(4))$ are isomorphic under the map
$$(f(z),h(z))\mapsto
f(4z)H_{\frac{7}{2}}(\tau)+h(4z)H_{\frac{11}{2}}(z).$$ For
$k\geq2$ we have
$M^{+}_{k+\frac{1}{2}}(\Gamma_0(4))=\mathbb{C}H_{k+\frac{1}{2}}\bigoplus
S^+_{k+\frac{1}{2}}(\Gamma_0(4))$.
\end{prop}

The results of
 \cite{KohZag}, \cite{Koh}, and \cite{wal} connect the
coefficients of Hecke eigenforms of half-integral weight to the
central $L$-values of twists of integral weight Hecke eigenforms.
More precisely, suppose that $f(z) = \sum_{n=1}^{\infty} a(n)
q^{n} \in S_{2k}$ is a normalized Hecke eigenform  and that $g(z)
= \sum_{n=1}^{\infty} c(n) q^{n} \in
S_{k+\frac{1}{2}}^{+}(\Gamma_{0}(4))$ is a Hecke eigenform with
the same Hecke eigenvalues as those of $f$. Here,

$$S_{k+\frac{1}{2}}^{+}(\Gamma_{0}(4)):=\{ g\in
S_{k+\frac{1}{2}}(\Gamma_0(4)) \, | \, g(z)=\sum_{n\geq 1,
(-1)^k n\equiv 0,1 \pmod{4}} c(n) q^n \}$$ Theorem~1 of
\cite{KohZag} states the following.

\begin{thm}\label{KohnenZagier}
Suppose that $f$ and $g$ are as above, $D$ is a fundamental
discriminant with $(-1)^{k} D > 0$, and $L(f,D,s)$ is the twisted
$L$-series
\[
  L(f,D,s) = \sum_{n=1}^{\infty} \left(\tfrac D n\right) a(n) n^{-s}.
\]
Then
\[
  \frac{c(|D|)^{2}}{\langle g, g \rangle} = \frac{(k-1)!}{\pi^{k}}
  |D|^{k - \frac{1}{2}} \frac{L(f,D,k)}{\langle f ,f \rangle}.
\]
\end{thm}
Here, $\langle g, g \rangle$ and $\langle f, f \rangle$ are the
normalized Petersson scalar products
\begin{align*}
  \langle g, g \rangle &= \frac{1}{6} \int_{\Gamma_{0}(4) \backslash \mathbb{H}}
  |g(z)|^{2} y^{k - 3/2} \, dx \, dy\\
  \langle f, f \rangle &= \int_{\Gamma_1 \backslash\mathbb{H} }
  |f(z)|^{2} y^{2k-2} \, dx \, dy.
\end{align*}

\bigskip

\section{ \bf{Proof of Theorem \ref{main}}}

 For any function $f(z)$ on $\mathbb H = \{ z \in
\C : {\rm Im}(z) > 0 \}$ and any positive integer $d$ we define
the operator $U_d$
\begin{gather}
(f \big|U_d)(z):=\frac1d\sum_{j=0}^{d-1}
f\left(\frac{z+j}d\right).\label{ud}
\end{gather}
Suppose that $f$ has a Fourier expansion $f(z)=\sum a(n)q^n$. Then
$(f\big|U_d)(z)=\sum a(nd)q^n$. If
$g(z)=\sum_{n=0}^{\infty}c(n)q^n$ is in
$M_{k+\frac{1}{2}}(\Gamma_0(4N)),$  then the Hecke operator with
the trivial character on $M_{k+\frac{1}{2}}(\Gamma_0(4N))$   is
defined for odd primes $\ell $ by
\begin{equation*}
(g  \big|
\,T(\ell^2,k))(z):=\sum_{n=0}^{\infty}\left(c(\ell^2n)+\ell^{k-1}\left(\frac{(-1)^kn}{\ell}\right)
c(n)+\left(\frac{(-1)^k}{\ell^2}\right)\ell^{2k-1}c\left(\frac{n}{\ell^2}\right)\right)q^n,
\end{equation*}
where $\left(\frac{\cdot}{\ell}\right)$ and
$\left(\frac{\cdot}{\ell^2}\right)$ are Jacobi symbols, and
$c\left(\frac{n}{\ell^2}\right):=0$ if $\ell^2\nmid n$. If $g$ has
integral coefficients, then one also has that
\begin{equation}
\label{U-T} g\, \big| \,U_{{\ell}}\equiv g^{\ell}\, \big|
\,T\left(\ell^2,\ell k+\frac{(\ell-1)}{2}\right)\pmod{\ell}.
\end{equation}

 For any positive integer $d$ we define the
operator $V_d$
\begin{gather}
(g \big|V_d)(z):=\sum_{n=0}^{\infty}c(n)q^{dn}.
\end{gather}
Note that if $\ell$ is a prime, then
\begin{gather}
(g \big|V_{\ell})(z)\equiv g(z)^{\ell} \pmod{\ell}.
\end{gather}
For a Dirichlet character $\chi$ let
$$g\otimes\chi:=\sum_{n=0}^{\infty}\chi(n)c(n)q^n.$$ The following
proposition immediately implies our main theorems.
\begin{prop}\label{mainprop}
Suppose that $k\geq 8$ is an integer.
\begin{enumerate}
\item If $k$ is odd, then for any positive $\epsilon$
\begin{equation*}
\mathcal{N}_{k,\Gamma_1}(x)\geq \left(\frac{9}{16\pi^2}-\epsilon\right)x \,
\, \,  (x\gg_{\epsilon} 0).
\end{equation*}
\item If $k$ is an even integer such that $d_k>1$ or $k=10$, then
for any positive $\epsilon$
\begin{equation*}
\mathcal{N}_{k,\Gamma_1}(x)\geq \left(\frac{9}{16\pi^2}-\epsilon\right)x \,
\, \,  (x\gg_{\epsilon} 0).
\end{equation*}
\end{enumerate}
\end{prop}

{\pf } For a positive even positive integer $t$, let
$$m(t):=\frac{1}{2}\left(t-4\left[\frac{t}{4}\right]\right).$$
For a non negative integer $t$ we define a  modular form $R_t(z)$:
$$
    R_t(z):=\left\{
    \begin{array}{cc}
    E_4(4z)^{\left[\frac{t}{4}\right]-m(t)}E_6(4z)^{m(t)}& \text{ if } t>0, \\
     1 & \text{ if } t=0.\end{array}\right.$$
For any even positive integer $t,$ we have
\begin{equation}\label{R}
R_t(z)\equiv 1\pmod{3}.
\end{equation}
This is from the fact that if $k\geq4$ is an
even integer, then $E_{p-1}(z) \equiv 1 \pmod{p}$ for any prime
$p$ such that $k\equiv0 \pmod{p-1}$ (see \cite{L}).

First, we assume that $k$ is odd. For each odd integers $k \geq 9$   define
\begin{equation}
\Phi_{k+\frac{1}{2}}(z)=28H_{3+\frac{1}{2}}(z)R_{k-3}(z)-\frac{44}{3}H_{5+\frac{1}{2}}(z)R_{k-5}(z)
:=\sum_{n=1}^{\infty}\beta_{k}(n)q^n.
\end{equation}
Then $\Phi_{k+\frac{1}{2}}(z)$ is in
$S_{k+\frac{1}{2}}^{+}(\Gamma_0(4))$ by Proposition
\ref{isomorphism}, and the Fourier coefficients of $\Phi_{k+\frac{1}{2}}(z)$
are 3-integral. On the other hand, for every odd $k\geq 9,$ we have
$$\Phi_{k+\frac{1}{2}}(z)\equiv\Phi_{9+\frac{1}{2}}(z)(=
 28H_{3+\frac{1}{2}}(z)R_{6}(z)-\frac{44}{3}H_{5+\frac{1}{2}}(z)R_{4}(z))
\pmod{3}.$$ Let
\begin{align*}
F(z)&:=(\Phi_{9+\frac{1}{2}}(z)-(\Phi_{9+\frac{1}{2}}(z)|U_3|V_3))
+(\Phi_{9+\frac{1}{2}}(z)-(\Phi_{9+\frac{1}{2}}(z)|U_3|V_3))
\otimes\left(\frac{\cdot}{3}\right)\\
&=\sum_{
\sm   \\
n \equiv 1 \pmod{3} \esm } \beta_9(n)q^n.
\end{align*}

Recall that $$G_{3,1}(z)=\sum_{ \sm n \equiv 1 \pmod{3} \esm }
H(1,n)q^n.$$

Proposition \ref{eis} implies that $G_{3,1}(z) =\sum_{ \sm n \equiv
1 \pmod{3} \esm } H(1,n)q^n $ is  a modular form of weight
$\frac{3}{2}$ on $\Gamma_0(36)$ such that its coefficients are
3-integral.

\medskip

Thus, by computing a few coefficient modulo $3$  and using Sturm's bound in (\cite{st}) we have
\begin{equation}\label{cong}
\begin{aligned}
F(z)&\equiv\sum_{  n \equiv 1 \pmod{3}   } \beta_9(n)q^n \\
&\equiv  2q^4+q^7+q^{19}+2q^{28}+2q^{40}+q^{43}+2q^{52}+q^{55}+2q^{64}+q^{67}+q^{76}
+\cdots\\
&\equiv \sum_{ n
\equiv 1 \pmod{3}   } H(1,n)q^n \equiv G_{3,1}(z) \pmod{3}.
\end{aligned}
\end{equation}
On the other hand, let $h(D)=H(1,D)$ be the class number of
$\mathbb{Q}(\sqrt{D}).$ It is known that for $D<0$,
$h(D)=-B_{1,\left(\frac{D}{\bullet}\right)}$ (apart from $D=-3$
and $-4$) (for example, see \cite{Urba}). Thus, we have
$$\beta_k(D)\equiv B_{1,(\frac{-D}{\cdot})} \equiv
h(-D) \pmod{3},$$ for a fundamental discriminant $ D>1$ such that
$D \equiv 1 \pmod{3}.$\\

Now let  $m$ and $N$ be positive integers satisfying the
condition: \\

\medskip
 (**) \, \, {\em{
 If an odd prime  $p$ is a common divisor of $m$  and $ N,$
then $p \mid N$ and $p^2 \nmid m$. Further if $N$ is even, then
$(i)$ $4|N$ and $m\equiv 1 \mod{4}$ or $(ii)$  $16|N$ and $m\equiv
8 \text{ or } 12 \mod{16}$.
}}
\medskip

We denote by $N_2^- (x,m,N)$ the number of fundamental
discriminants $D$ with $-x<D<0$ and $D \equiv m \pmod{N}$. The
results of Davenport-Heilbronn \cite{D-H} and Nakagawa-Horie
\cite{Na-Ho} imply that for any positive number $\epsilon$
\begin{equation}\label{class number}
\begin{array}{ll}
 |\{  \text{ fundamental discriminants } D\equiv 1 \pmod{3} \; | \; 0 < D < x
\text{ and } 3 \nmid h(-D) \}|\\
 \gg \left(\frac{1}{2}-\epsilon \right)(N_2^-(x,1,3)).
\end{array}
\end{equation}
Since $N_2^-(x,1,3)\sim \frac{9}{8\pi^2}x$ for $x\rightarrow
\infty$ (see Proposition 2. in \cite{Na-Ho}), for odd integers
$k\geq 7$ we have
\begin{equation}\label{nonvanishing}
|\{ \text{ fundamental discriminants } D \; | \; 0 < D < x \text{
and } 3 \nmid \beta_k(D)\}| \gg \left(\frac{9}{16 \pi^2}-\epsilon
\right)x.
\end{equation}
Since $\Phi_{k+\frac{1}{2}}(z) \in
S_{k+\frac{1}{2}}^{+}(\Gamma_{0}(4))$, the cusp form
$\Phi_{k+\frac{1}{2}}(z)$ is a linear combination of Hecke
eigenforms $g_{\ell}(z)=\sum_{n=1}^{\infty}c_{\ell}(n)q^n \in
S_{k+\frac{1}{2}}^{+}(\Gamma_{0}(4))$ for $1 \leq  \ell \leq d_k$.
By Theorem\ref{KohnenZagier} we complete the proof of Theorem.
\\

From now on, suppose that $k$ is even, and that $d_k>1$ or $k=10$. Let
\begin{eqnarray}{\label{9}}
\Psi_{k+\frac{1}{2}}(z):=\Delta(4z)R_{k-12}(z)\theta(z)
=\sum_{n=1}^{\infty} \alpha_{k}(n)q^n \;\; \text{ for } k>10
\end{eqnarray}
 and
\begin{eqnarray}{\label{10}}
\Psi_{10+\frac{1}{2}}(z):= -(\theta(z)E_4(4z)E_6(4z)
-H_{\frac{5}{2}}(z)E_4(4z)^2)\otimes\chi_3\end{eqnarray}
$$
+(\theta(z)E_4(4z)E_6(4z)-H_{\frac{5}{2}}(z)E_4(4z)^2)\otimes\chi_3^2,$$
where $\chi_3(n)=\left(\frac{n}{3}\right)$. We have by
the sturm's bound and (\ref{R})

$$
\sum_{n \equiv 2 \pmod{3}} \alpha_{k}(n)q^n  \equiv \sum_{n \equiv
2 \pmod{3}} \alpha_{12}(n)q^n$$ $$\equiv 2 q^8+2q^{17}+ q^{20} +
2q^{41} + q^{44}+q^{53}+q^{56}+q^{65}+2q^{68}+2q^{80} +2q^{89} +
2q^{92} \cdots $$ $$ \equiv  \sum_{n \equiv 2 \pmod{3}} H(1,3n)q^n
\pmod{3}.
$$
The remained part of the proof can be completed in a  similar way
as before, so we omit the details. \qed

\begin{rem} The argument given in  [\cite{Kohnen3}, p. $186$
bottom] in the case where $k$ is even and $k\equiv 1 \pmod 3$ is
not correct, since it would require that all the coefficients of
$\delta_{k-4}$ are 3-integral which in general is not the case.
\end{rem}

\bigskip


\section{\bf{Conclusion}}


In this paper, we extend the result in \cite{Kohnen3} to the case
when $k$ is any integer, in particular when $k$ is odd. So, we
obtain that, for each integer $2k$ such that dim$S_{2k}\geq 1,$  there
is a normalized Hecke eigenform $F$ in $S_{2k}$ satisfying $
\mathcal{N}_F(x)\gg x \,\, \, \, (x \rightarrow \infty). $ We
conclude this paper with the following remark:

\begin{rem}
\begin{enumerate}
\item For each odd $k\geq9$ and even $\lambda$ such that
$\lambda\geq 12$, all the coefficients of
$\Phi_{k+\frac{1}{2}}$ and $\Psi_{\lambda+\frac{1}{2}}$ are
3-integral and a positive portion of these coefficients
$\beta_k(n)$ and $\alpha_{\lambda}(n)$ is not vanishing modulo
$3$.

\item Note that $k=9$ is the minimum odd integer such that
$\dim(S_{2k})>0$. Let
$$\Delta(z)=q\prod_{n=1}^{\infty}(1-q^n)^{24}.$$ If we take
$F(z)=\Delta(z)E_6(z)$, then $F(z)$ is the unique normalized Hecke
eigenform in $S_{18}.$ Corollary\ref{sp} implies that
\begin{equation*}
\mathcal{N}_F(x)\geq \left(\frac{9}{16\pi^2}-\epsilon\right)x \,
\, \,  (x\gg_{\epsilon} 0).
\end{equation*}

\item The direct computation shows that if
 $f\in
S_{k+\frac{1}{2}}^+(\Gamma_0(4))$ has integral coefficients for
even $k$ such that $d_{k}=1$, then
$$f \equiv c \sum_{\sm n\geq 1 \\ 3 \nmid n \esm}q^{n^2} \pmod{3}$$
for some $c  \in \{-1,1\}$.

\end{enumerate}
\end{rem}

\section*{\bf{Acknowledgement}}
The first author was supported by   KRF-2008-331-C00005  and wish
to express his gratitude to KIAS for its support through Associate
membership program. The second author was partially supported by
KOSEF-R01-2008-000-20446 and KRFKRF-2007-412-J02302.

\end{document}